\documentclass[12pt,a4paper]{article}
\usepackage[english]{babel}
\usepackage{amsmath}
\usepackage{amssymb}
\usepackage{amsfonts}
\usepackage{hyperref}

\begin{document}

\binoppenalty=10000 \relpenalty=10000

\renewcommand{\refname}{References}
\renewcommand{\contentsname}{Contents}

\begin{center}
{\huge Global solvability of the initial boundary value problem for a model system of one-dimensional equations of polytropic flows of viscous compressible fluid mixtures}
\end{center}

\medskip

\begin{center}
{\large Dmitriy Prokudin\footnote{This research was supported by the Ministry of Education and Science of the Russian Federation (grant 14.Z50.31.0037).}}
\end{center}

\medskip

\begin{center}
{\large October 20, 2017}
\end{center}

\medskip

\begin{center}
{
Lavrentyev Institute of Hydrodynamics, 630090 Novosibirsk, Russia\\
Voronezh State University, 394018 Voronezh, Russia}
\end{center}

\medskip

\begin{center}
{\bfseries Abstract}
\end{center}


\begin{center}
\begin{minipage}{110mm}
We consider the initial boundary value problem for a model system of one-dimensional equations which describe unsteady polytropic motions of a mixture of viscous compressible fluids. We prove the global existence and uniqueness theorem for the strong solution without restrictions on the structure of the viscosity matrix except standard properties of symmetry and positiveness.
\end{minipage}
\end{center}

\bigskip


\newpage

\tableofcontents

\bigskip

\section{Statement of the problem, formulation\\
of the result, the Lagrangian coordinates}

\noindent\indent We consider the problem of one-dimensional polytropic flow of a mixture which consists of $N\geqslant 2$ components. We work in the closure $\overline{Q}_{T}$ of the domain $Q_{T}=(0, 1)\times (0, T)$, where $T>0$~is an arbitrary positive number, and our aim is to find the density $\rho> 0$ of the mixture and the velocity $u_{i}$ for each component of the mixture numbered by $i=1,\ldots,N$, which satisfy the following system of equations, initial and boundary conditions \cite{prok17.mamprok17}:
\begin{equation}\label{prok17.continuity}\partial_{t}\rho+\partial_{x}(\rho v)=0,\quad v=\frac{1}{N}\sum\limits_{i=1}^{N}u_{i},\end{equation}
\begin{equation}\label{prok17.momentum}
\rho\left(\partial_{t}u_{i}+v\partial_{x}u_{i}\right)+K\partial_{x} \rho^{\gamma}=\sum\limits_{j=1}^N \mu_{ij}\partial_{xx}u_{j}+\sum\limits_{j=1}^{N}a_{ij}(u_{j}-u_{i}),\quad i=1,\ldots,N,
\end{equation}
\begin{equation}\label{prok17.nachusl}\rho|_{t=0}=\rho_{0}, \quad u_{i}|_{t=0}=u_{0i},\quad i=1,\ldots,N,\end{equation}
\begin{equation}\label{prok17.boundvelocity}u_{i}|_{x=0}=u_{i}|_{x=1}=0,\quad i=1,\ldots,N.\end{equation}
Here $v$ is the average velocity of the mixture; the values $K>0$, $\displaystyle \gamma>1$, $\mu_{ij}=\mu_{ji}$ and $a_{ij}=a_{ji}>0$, $i, j =1,\ldots,N$ are known constants, and the viscosity coefficients $\{\mu_{ij}\}_{i, j=1}^{N}$ compose the matrix $\textbf{M}>0$; the initial distributions $\rho_{0}$ and $u_{0i}$, $i=1,\ldots, N$, are prescribed.

The aim of the paper is to prove the existence and uniqueness of the strong solution to the problem (\ref{prok17.continuity})--(\ref{prok17.boundvelocity}).

{\bfseries Definition 1.} {\it Strong solution to the problem $($\ref{prok17.continuity}$)$--$($\ref{prok17.boundvelocity}$)$ is called the collection of $N+1$ functions $(\rho, u_{1},\ldots, u_{N})$
such that the equations $($\ref{prok17.continuity}$)$, $($\ref{prok17.momentum}$)$ are satisfied almost everywhere in $Q_{T}$, the initial data $($\ref{prok17.nachusl}$)$ are accepted for a.a. $x\in (0, 1)$, the boundary conditions $($\ref{prok17.boundvelocity}$)$ are valid for a.a. $ t\in (0, T)$, and the following inequality and inclusions hold $(i=1,\ldots,N)$
\begin{equation}\label{prok17.oprresh}
\begin{array}{c}
\displaystyle
\rho>0,\quad \rho\in L_{\infty}\big(0, T; W^{1}_{2}(0, 1)\big), \quad  \partial_{t}\rho\in L_{\infty}\big(0, T; L_{2}(0, 1)\big),\\ \\
\displaystyle u_{i}\in L_{\infty}\big(0, T; W^{1}_{2}(0, 1)\big)\bigcap L_{2}\big(0, T; W^{2}_{2}(0, 1)\big),\quad
\partial_{t}u_{i} \in L_{2}(Q_{T}).
\end{array}
\end{equation}}

The main result of the paper is formulated as the following theorem.

{\bfseries Theorem 2.} {\it Let the initial data in $($\ref{prok17.nachusl}$)$ satisfy the conditions\linebreak  $(i=1,\ldots,N)$
\begin{equation}\label{prok17.gladknachusl}
\rho_{0}\in W^{1}_{2}(0, 1),\; \rho_{0}>0,\; u_{0i}\in {W^1_2}(0, 1),\; u_{0i}|_{x=0}=u_{0i}|_{x=1}=0,
\end{equation}
symmetric viscosity matrix $\textbf{M}$ is positive, the polytropic exponent $\gamma>1$, all other numeric parameters $K$, $T$ and $a_{ij}=a_{ji}$, $i, j =1,\ldots,N$, are positive.

Then there exists the unique strong solution to the problem $($\ref{prok17.continuity}$)$--$($\ref{prok17.boundvelocity}$)$ in the sense of Definition~1.}

{\bf Sketch of the proof of Theorem 2.} The existence of unique strong solution to the problem (\ref{prok17.continuity})--(\ref{prok17.boundvelocity}) in a small time interval $[0, t_{0}]$ is proved in \cite{prok17.prok17}. In order to extend this solution from the interval $[0, t_{0}]$ to the target interval $[0, T]$, we need to prove a priori estimates, in which the constants depend only on the input data of the problem and on the value $T$, but not on the small parameter $t_{0}$ (see, e.~g., \cite{prok17.akm90}). That is why we concentrate on the global estimates.

{\bf Lagrangian coordinates.} While the problem (\ref{prok17.continuity})--(\ref{prok17.boundvelocity}) is studied, it is sometimes more convenient to use the Lagrangian coordinates. Let us consider  $\displaystyle y(x,t)=\int\limits_{0}^{x}\rho(s,t)\,ds$ and $t$ as new independent variables. Then the system (\ref{prok17.continuity}), (\ref{prok17.momentum}) turns into the form
\begin{equation}\label{prok17.1newcontinuity1lagr}
\partial_{t}\rho+\rho^{2}\partial_{y}v=0,
\end{equation}
\begin{equation}\label{prok17.1newmomentum1lagr}
\partial_{t}u_{i}+K\partial_{y} \rho^{\gamma}=
\sum\limits_{j=1}^N \mu_{ij}\partial_{y}(\rho\partial_{y}u_{j})+\frac{1}{\rho}\sum\limits_{j=1}^{N}a_{ij}(u_{j}-u_{i}),\quad i=1,\ldots,N.
\end{equation}
The domain $Q_{T}$ is mapped into the rectangular $\Pi_{T}=(0, d)\times(0, T)$, where $\displaystyle d=\int\limits_{0}^{1}\rho_{0}(x)\,dx>0$, and the initial and boundary conditions accept the form
\begin{equation}\label{prok17.1nachusl1lagr}
\rho|_{t=0}=\widetilde{\rho}_{0}, \quad u_{i}|_{t=0}=\widetilde{u}_{0i},\quad i=1,\ldots,N,
\end{equation}
\begin{equation}\label{prok17.1boundvelocity1lagr}
u_{i}|_{y=0}=u_{i}|_{y=d}=0,\quad i=1,\ldots,N.
\end{equation}

\section{Global a priori estimates}

\noindent\indent Let us multiply the equations (\ref{prok17.momentum}) by $ u_{i}$, integrate the result over $(0, 1)$ and sum over $i=1,\ldots,N$. Due to (\ref{prok17.continuity}), (\ref{prok17.boundvelocity}) and the condition $\textbf{M}>0$, the following relations hold
\begin{equation}\label{prok17.lemma51.2}
\sum\limits_{i=1}^{N}\int\limits^{1}_{0}\Big(\rho\partial_{t}u_{i}+\rho v\partial_{x}u_{i}\Big)u_{i}\,
dx= \frac{1}{2}\frac{d}{dt}\left(\sum\limits_{i=1}^{N}
\int\limits^{1}_{0}\rho u_{i}^{2}\, dx\right),
\end{equation}
\begin{equation}\label{prok17.lemma51.3}
K\sum\limits_{i=1}^{N}\int\limits^{1}_{0}u_{i}\left(\partial_{x}\rho^{\gamma}\right)\,
dx=-KN\int\limits^{1}_{0}\rho^{\gamma}\left(\partial_{x}v\right)\, dx=
\frac{KN}{\gamma-1}\frac{d}{dt}\left(\int\limits^{1}_{0}\rho^{\gamma}\, dx\right),
\end{equation}
\begin{equation}\label{prok17.lemma51.4}
\begin{array}{c}
\displaystyle
\sum\limits_{i,j=1}^N
\mu_{ij}\int\limits\limits_{0}^{1}(\partial_{xx}u_{j})u_{i}\,
dx=-\sum\limits_{i,j=1}^N
\mu_{ij}\int\limits^{1}_{0}(\partial_{x}u_{i})(\partial_{x}u_{j})\,
dx\leqslant\\ \\
\displaystyle
\leqslant -C_{0}(\textbf{M})
\sum\limits_{i=1}^N\int\limits\limits_{0}^{1}|\partial_{x}u_{i}|^{2}\,
dx,
\end{array}
\end{equation}
\begin{equation}\label{prok17.newlemmanew51.5}
\sum\limits_{i, j=1}^N
a_{ij}\int\limits\limits_{0}^{1}(u_{j}-u_{i})u_{i}\,
dx=-\frac{1}{2}\sum\limits_{i, j=1}^{N} a_{ij}\int\limits_{0}^{1}(u_{i}-u_{j})^{2}\, dx,
\end{equation}
and we come to the inequality
\begin{equation}\label{prok17.lemma51.6}
\begin{array}{c}
\displaystyle
\frac{d}{dt}\sum\limits_{i=1}^{N}
\int\limits^{1}_{0}\left(\frac{1}{2}\rho u_{i}^{2}+\frac{K}{\gamma-1}\rho^{\gamma}\right)\,
dx+
C_{0}\sum\limits_{i=1}^N\int\limits\limits_{0}^{1}|\partial_{x}u_{i}|^{2}\,
dx+\\ \\
\displaystyle
+\frac{1}{2}\sum\limits_{i, j=1}^{N} a_{ij}\int\limits_{0}^{1}(u_{i}-u_{j})^{2}\, dx\leqslant 0.
\end{array}
\end{equation}
Let us agree that $C_k(\cdot)$, $k=0, 1, \ldots$, stand for the quantities which take finite positive values and depend on the objects listed in the parentheses. When we integrate the inequality (\ref{prok17.lemma51.6}) over $(0, t)$ using (\ref{prok17.nachusl}), we obtain the estimate
\begin{equation}\label{prok17.lemma51.7}
\begin{array}{c}
\displaystyle \sum\limits_{i=1}^{N}
\int\limits^{1}_{0}\left(\frac{1}{2}\rho u_{i}^{2}+\frac{K}{\gamma-1}\rho^{\gamma}\right)\,
dx+
C_{0}\sum\limits_{i=1}^N\int\limits\limits_{0}^{t}\int\limits\limits_{0}^{1}|\partial_{x}u_{i}|^{2}\,
dxd\tau+ \\ \\
\displaystyle +\frac{1}{2}\sum\limits_{i, j=1}^{N} a_{ij}\int\limits\limits_{0}^{t}\int\limits_{0}^{1}(u_{i}-u_{j})^{2}\, dxd\tau\leqslant \sum\limits_{i=1}^{N}
\int\limits^{1}_{0}\left(\frac{1}{2}\rho_{0}u_{0i}^{2}+\frac{K}{\gamma-1}\rho_{0}^{\gamma}\right)\,
dx,
\end{array}
\end{equation}
which leads to the fact that
\begin{equation}\label{prok17.lemma1}
\begin{array}{c}
\displaystyle\sum\limits_{i=1}^{N}\left(\|\sqrt{\rho}u_{i}\|_{L_{\infty}\big(0, T;
L_{2}(0, 1)\big)}
+\|\partial_{x}u_{i}\|_{L_{2}(Q_{T})}+\sum\limits_{j=1}^{N}\|u_{i}-u_{j}\|_{L_{2}(Q_{T})}\right)+\\ \\
\displaystyle
+\|\rho\|_{L_{\infty}\big(0, T;
L_{\gamma}(0,
1)\big)}\leqslant
C_{1},
\end{array}
\end{equation}
where $C_{1}=C_{1}\left(\left\{\|\sqrt{\rho_{0}}u_{0i}\|_{L_{2}(0,
1)}\right\}, \|\rho_{0}\|_{L_{\gamma}(0,
1)}, \{a_{ij}\}, K, \textbf{M}, N, \gamma\right)$. The estimates (\ref{prok17.lemma1}), due to (\ref{prok17.1boundvelocity1lagr}), lead to the inequality
\begin{equation}\label{prok17.0806171}
\sum\limits_{i=1}^{N}\|u_{i}\|_{L_{2}\big(0, T;
L_{\infty}(0, 1)\big)}\leqslant C_{1}.
\end{equation}
The estimate (\ref{prok17.lemma1}) in the Lagrangian coordinates takes the form
\begin{equation}\label{prok17.lemma1lagr}
\begin{array}{c}
\displaystyle
\sum\limits_{i=1}^{N}\Big(\|u_{i}\|_{L_{\infty}\big(0, T;
L_{2}(0, d)\big)}
+\|\sqrt{\rho}\partial_{y}u_{i}\|_{L_{2}(\Pi_{T})}+\\ \\
\displaystyle +\sum\limits_{j=1}^{N}\|(u_{i}-u_{j})/\sqrt{\rho}\|_{L_{2}(\Pi_{T})}\Big)+\|\rho\|_{L_{\infty}\big(0, T;
L_{\gamma-1}(0,
d)\big)}\leqslant
C_{2}(C_{1}, \gamma).
\end{array}
\end{equation}

The next step is to prove the positiveness and boundedness of the density~$\rho$. Here we use the equations (\ref{prok17.1newcontinuity1lagr}), (\ref{prok17.1newmomentum1lagr}). Let us rewrite the equations (\ref{prok17.1newmomentum1lagr}) in the form
\begin{equation}\label{prok17.newnewmomentum}
\begin{array}{c}
\displaystyle
\sum\limits_{j=1}^{N}\widetilde{\mu}_{ij}\partial_{t}u_{j}+K\left(\sum\limits_{j=1}^{N}\widetilde{\mu}_{ij}\right)\partial_{y} \rho^{\gamma}=
\partial_{y}(\rho\partial_{y}u_{i})+\\ \\
\displaystyle
+\frac{1}{\rho}\sum\limits_{j=1}^{N}\widetilde{\mu}_{ij}\left(\sum\limits_{k=1}^{N}a_{jk}(u_{k}-u_{j})\right),\quad i=1,\ldots,N,
\end{array}
\end{equation}
where $\widetilde{\mu}_{ij}$ are the entries of the matrix $\widetilde{\textbf{M}}=\textbf{M}^{-1}>0$, then sum (\ref{prok17.newnewmomentum}) over $i=1,\ldots,N$, and divide by $N$. Then we come to the equality
\begin{equation}\label{prok17.newnewmomentum1}
\partial_{t} V+\widetilde{K}\partial_{y}\rho^{\gamma}=
\partial_{y}\left(\rho\partial_{y}v\right)+\frac{1}{N\rho}\sum\limits_{i, j=1}^{N}\widetilde{\mu}_{ij}\left(\sum\limits_{k=1}^{N}a_{jk}(u_{k}-u_{j})\right),
\end{equation}
where $\displaystyle V=\frac{1}{N}\sum\limits_{i, j=1}^{N}\widetilde{\mu}_{ij}u_{j}$ and $\displaystyle \widetilde{K}=\frac{K}{N}\sum\limits_{i, j=1}^{N}\widetilde{\mu}_{ij}>0$.
Let us use (\ref{prok17.1newcontinuity1lagr}) in order to express
\begin{equation}\label{prok17.eq0102173}
\rho\partial_{y}v=-\partial_{t}\ln\rho,
\end{equation}
and substitute this relation into (\ref{prok17.newnewmomentum1}):
\begin{equation}\label{prok17.eq0102171}
\partial_{ty}\ln\rho+\widetilde{K}\partial_{y}\rho^{\gamma}=-\partial_{t} V
+\frac{1}{N\rho}\sum\limits_{i, j=1}^{N}\widetilde{\mu}_{ij}\left(\sum\limits_{k=1}^{N}a_{jk}(u_{k}-u_{j})\right).
\end{equation}
Let us multiply this equality by $\displaystyle \partial_{y}\ln\rho=:w$ and integrate over $y\in(0,d)$, then we obtain the equality
\begin{equation}\label{prok17.eq0102172}
\begin{array}{c}
\displaystyle
\frac{1}{2}\frac{d}{dt}\left(\int\limits_{0}^{d}w^{2}\, dy\right)+\widetilde{K}\gamma\int\limits_{0}^{d}\rho^{\gamma}w^{2}\,dy
=-\int\limits_{0}^{d}\left(\partial_{t} V\right)w\,dy
+\\ \\
\displaystyle
+\frac{1}{N}\sum\limits_{i, j, k=1}^{N}\widetilde{\mu}_{ij}a_{jk}\int\limits_{0}^{d}\frac{(u_{k}-u_{j})w}{\rho}\,dy.
\end{array}
\end{equation}
We transform the first summand in the right-hand side of (\ref{prok17.eq0102172}) via integration by parts and using (\ref{prok17.eq0102173}):
\begin{equation}\label{prok17.eq01021777}
-\int\limits_{0}^{d}\left(\partial_{t} V\right)w\,dy
=-\frac{d}{dt}\left(\int\limits_{0}^{d} Vw\,dy\right)+\int\limits_{0}^{d} \rho(\partial_{y}v)(\partial_{y}V)\,dy,
\end{equation}
and we estimate the second summand from above:
\begin{equation}\label{prok17.eq0102174}
\begin{array}{c}
\displaystyle
\frac{1}{N}\sum\limits_{i, j, k=1}^{N}\widetilde{\mu}_{ij}a_{jk}\int\limits_{0}^{d}\frac{(u_{k}-u_{j})w}{\rho}\,dy\leqslant
C_{3}\|1/\sqrt{\rho}\|_{L_{\infty}(0,d)}\|w\|_{L_{2}(0,d)}\times \\
\displaystyle
\times\sum\limits_{j, k=1}^{N}\|(u_{k}-u_{j})/\sqrt{\rho}\|_{L_{2}(0,d)},
\end{array}
\end{equation}
where $C_{3}=C_{3}\left(\{a_{jk}\}, \widetilde{\textbf{M}}\right)$. It is obvious, due to the equation (\ref{prok17.1newcontinuity1lagr}) and the conditions (\ref{prok17.1nachusl1lagr}) and (\ref{prok17.1boundvelocity1lagr}), that for every $t\in[0,T]$
\begin{equation}\label{prok17.eq0102176}
\rho(z(t), t)=d
\end{equation}
at least in one point $z(t)\in[0, d]$. Hence, we can use the representation
\begin{equation}\label{prok17.eq0102177}
\begin{array}{c}
\displaystyle
\frac{1}{\sqrt{\rho(y,t)}}=\frac{1}{\sqrt{\rho(z(t),t)}}+\int\limits_{z(t)}^{y}\partial_{s}\rho^{-\frac{1}{2}}(s,t)\, ds=\\ \\
\displaystyle
=d^{-\frac{1}{2}}-\frac{1}{2}\int\limits_{z(t)}^{y}\rho^{-\frac{1}{2}}(s,t)
\partial_{s}\ln{\rho(s,t)}\, ds,
\end{array}
\end{equation}
from which, using H\"older's inequality and (\ref{prok17.eq0102176}), we obtain
\begin{equation}\label{prok17.eq0102178}
\|1/\sqrt{\rho}\|_{L_{\infty}(0,d)}\leqslant d^{-\frac{1}{2}}+\frac{1}{2}\|w\|_{L_{2}(0,d)}.
\end{equation}
Hence, after the integration of (\ref{prok17.eq0102172}) over $(0,t)$, using (\ref{prok17.eq01021777}), (\ref{prok17.eq0102174}) and (\ref{prok17.eq0102178}), we come to the inequality
\begin{equation}\label{prok17.eq0102179}
\begin{array}{c}
\displaystyle \|w\|^{2}_{L_{2}(0, d)}+2\widetilde{K}\gamma\int\limits_{0}^{t}\int\limits_{0}^{d}\rho^{\gamma}w^{2}\,dyd\tau
\leqslant\|w_{0}\|^{2}_{L_{2}(0, d)}-2\int\limits_{0}^{d}Vw\,dy+\\\\
\displaystyle +2\int\limits_{0}^{d}V_{0}w_{0}\,dy+2\int\limits_{0}^{t}\|\rho(\partial_{y}V)(\partial_{y}v)\|_{L_{1}(0, d)}\,d\tau+\\\\
\displaystyle+C_{3}\sum\limits_{j,k=1}^{N}\int\limits_{0}^{t}\|(u_{k}-u_{j})/\sqrt{\rho}\|_{L_{2}(0, d)}\|w\|_{L_{2}(0,d)}\big(2d^{-\frac{1}{2}}+\|w\|_{L_{2}(0,d)}\big)\,d\tau,
\end{array}
\end{equation}
where $w_{0}=w(0,t)$ and $V_{0}=V(0,t)$. Using Cauchy's inequality and the estimate (\ref{prok17.lemma1lagr}), we derive from the last formula that
\begin{equation}\label{prok17.eq01021710}
\|w\|^{2}_{L_{2}(0, d)}\leqslant C_{4}+C_{5}\sum\limits_{j,k=1}^{N}\int\limits_{0}^{t}\|(u_{k}-u_{j})/\sqrt{\rho}\|_{L_{2}(0, d)}\|w\|^{2}_{L_{2}(0,d)}\,d\tau,
\end{equation}
where $C_{4}=C_{4}\left(C_{2}, C_{3}, \{\|\widetilde{u}_{0j}\|_{L_{2}(0,d)}\}, \|w_{0}\|_{L_{2}(0,d)}, \widetilde{\textbf{M}}, N, T, d\right)$ and $C_{5}=C_{5}(C_{3})$. Since (\ref{prok17.lemma1lagr}) leads to the estimate
\begin{equation}\label{prok17.eq01021711}
\sum\limits_{j,k=1}^{N}\int\limits_{0}^{t}\|(u_{k}-u_{j})/\sqrt{\rho}\|_{L_{2}(0, d)}\,d\tau\leqslant C_{6}(C_{2}, T)\quad \forall\, t\in[0,T],
\end{equation}
then the Gronwall lemma provides
\begin{equation}\label{prok17.eq01021712}
\|w(t)\|_{L_{2}(0,d)}\leqslant C_{6}(C_{4}, C_{5}, C_{6})\quad \forall\, t\in[0,T],
\end{equation}
i.~e. the norm of the derivative $\partial_{y}\ln\rho$ in $L_{2}(0,d)$ is bounded uniformly in $t\in[0,T]$. Hence, due to the representation (see the proof of (\ref{prok17.eq0102177}))
$$
{\ln{\rho(y,t)}}={\ln{\rho(z(t),t)}}+\int\limits_{z(t)}^{y}\partial_{s}\ln{\rho(s,t)}\, ds,
$$
we have
\begin{equation}\label{prok17.eq01021713}
|\ln\rho(y,t)|\leqslant |\ln{d}|+\sqrt{d}\|w\|_{L_{2}(0,d)}\leqslant C_{7}(C_{6},d),
\end{equation}
and consequently
\begin{equation}\label{prok17.eq01021714}
0<C_{8}^{-1}(C_{7})\leqslant\rho(y,t)\leqslant C_{8}(C_{7}).
\end{equation}

Now we possess the boundedness and positiveness of the density $\rho$, and the remaining a priori estimates can be obtained in the original (Eulerian) coordinates $(x,t)$. Thus, from (\ref{prok17.eq01021712}) and (\ref{prok17.eq01021714}) we deduce
\begin{equation}\label{prok17.eq01021715}
\left\|\partial_{x}\rho(t)\right\|_{L_{2}(0,1)}\leqslant C_{9}(C_{6}, C_{8})\quad \forall\, t\in[0,T].
\end{equation}
Then, let us square the momentum equations (\ref{prok17.momentum}) and sum over $i=1,\ldots,N$, as a result we obtain
\begin{equation}\label{prok17.1303171}\begin{array}{c}\displaystyle
\sum\limits_{i=1}^{N}\rho(\partial_{t}u_{i})^{2}+\frac{1}{\rho}\sum\limits_{i=1}^{N}\left(\sum\limits_{j=1}^{N}\mu_{ij}\partial_{xx}u_{j}\right)^{2}-
2\sum\limits_{i=1}^{N}(\partial_{t}u_{i})\left(\sum\limits_{j=1}^{N}\mu_{ij}\partial_{xx}u_{j}\right)=\\ \\
\displaystyle=\sum\limits_{i=1}^{N}\frac{1}{\rho}\left(\sum\limits_{j=1}^{N}a_{ij}(u_{j}-u_{i})-K\partial_{x}\rho^{\gamma}-\rho v\partial_{x}u_{i}\right)^{2}.
\end{array}
\end{equation}
Let us introduce the function $\alpha$ as
$$
\alpha(t)=\sum\limits_{i, j=1}^{N}\mu_{ij}\int\limits_{0}^{1}(\partial_{x}u_{i})(\partial_{x}u_{j})\, dx+$$$$+\sum\limits_{i=1}^{N}\int\limits_{0}^{t}\int\limits_{0}^{1}\left(\rho(\partial_{t}u_{i})^{2}+\frac{1}{\rho}\left(\sum\limits_{j=1}^{N}\mu_{ij}\partial_{xx}u_{j}\right)^{2}\right)\, dxd\tau.
$$
Then (\ref{prok17.1303171}) and the inequalities (\ref{prok17.lemma1}), (\ref{prok17.eq01021714}), (\ref{prok17.eq01021715}) lead to the estimate
\begin{equation}\label{prok17.1303172}
\begin{array}{c}
\displaystyle
\alpha^{\prime}(t)\leqslant C_{10}+C_{11}\left(\sum\limits_{j=1}^{N}\|u_{j}\|^{2}_{L_{\infty}(0,1)}\right)\left(\sum\limits_{i, j=1}^{N}\mu_{ij}\int\limits_{0}^{1}(\partial_{x}u_{i})(\partial_{x}u_{j})\, dx\right)\leqslant\\ \\
\displaystyle
\leqslant C_{10}+C_{11}\left(\sum\limits_{j=1}^{N}\|u_{j}\|^{2}_{L_{\infty}(0,1)}\right)\alpha(t),
\end{array}
\end{equation}
where $C_{10}=C_{10}(C_{1}, C_{8}, C_{9}, \{a_{ij}\}, K, N, \gamma)$ and $C_{11}=C_{11}(C_{8}, \textbf{M})$, from which via the Gronwall lemma (see also (\ref{prok17.0806171})) it follows that
$$\alpha(t)\leqslant C_{12}\left(C_{1}, C_{10}, C_{11}, \{\|u_{0i}^{\prime}\|_{L_{2}(0, 1)}\}, \textbf{M}, N, T\right).
$$
Using this and (\ref{prok17.eq01021714}), we come to the inequality
\begin{equation}\label{prok17.eq01021716}
\sum\limits_{i=1}^{N}\left(\|\partial_{x}u_{i}\|_{L_{\infty}(0, T; L_{2}(0,1))}+\|\partial_{xx}u_{i}\|_{L_{2}(Q_{T})}+\|\partial_{t}u_{i}\|_{L_{2}(Q_{T})}\right)\leqslant C_{13},
\end{equation}
where $C_{13}=C_{13}\left(C_{8}, C_{12}, \textbf{M}, N\right)$. Finally, the continuity equation (\ref{prok17.continuity}) and the inequalities (\ref{prok17.eq01021714}),  (\ref{prok17.eq01021715}) and  (\ref{prok17.eq01021716}) provide
\begin{equation}\label{prok17.eq01021717}
\|\partial_{t}\rho\|_{L_{\infty}(0, T; L_{2}(0,1))}\leqslant C_{14}(C_{8}, C_{9}, C_{13}).
\end{equation}
Theorem 2 is proved.

\end{document}